# Dynamic UAV-based traffic monitoring under uncertainty as a stochastic arc-inventory routing policy

Joseph Y.J. Chow
Department of Civil and Urban Engineering, New York University, New York, NY, USA

**Abstract**
Given the rapid advances in unmanned aerial vehicles, or drones, and increasing need to monitor at a city level, one of the current research gaps is how to systematically deploy drones over multiple periods. We propose a real-time data-driven approach: we formulate the first deterministic arc-inventory routing problem and derive its stochastic dynamic policy. The policy is expected to be of greatest value in scenarios where uncertainty is highest and costliest, such as city monitoring during major events. The Bellman equation for an approximation of the proposed inventory routing policy is formulated as a selective vehicle routing problem. We propose an approximate dynamic programming algorithm based on Least Squares Monte Carlo simulation to find that policy. The algorithm has been modified so that the least squares dependent variable is defined to be the "expected stock out cost upon the next replenishment". The new algorithm is tested on 30 simulated instances of real time trajectories over 5 time periods of the selective VRP to evaluate the proposed policy and algorithm. Computational results on the selected instances show that the algorithm on average outperforms the myopic policy by 23% to 28%, depending on the parametric design. Further tests are conducted on classic benchmark arc routing problem instances. The 11-link instance gdb19 is expanded into a sequential 15-period stochastic dynamic example and used to demonstrate why a naïve static multi-period deployment plan would not be effective in real networks.

# 1. Introduction

With the rise of Big Data analytics and urban informatics, there is an increasing interest to gather ever more real time data from a city's environment for real time traffic monitoring (Geroliminis and Daganzo, 2008), travel activity monitoring (Liu et al., 2014; Jiang et al., 2015), or humanitarian logistics (Ozguven and Ozbay, 2015), among others. Numerous monitoring sensor technologies exist for this purpose; some of the more promising among these are mobile sensors that can be deployed autonomously, such as unmanned aerial vehicles (UAVs, i.e. drones) (Chen et al., 2007). For example, UAVs have been demonstrated as feasible tools for gathering real traffic and transportation data (Srinivasan et al., 2004). UAVs can substitute traditional methods for a number of uses in transportation including measuring level of service, average annual daily traffic, intersection operations, parking utilization (Coifman et al., 2006); traffic management (Huiyuan et al., 2007); origin-destination estimation (Braut et al., 2012); and goods delivery as a "flying sidekick" (Murray and Chu, 2015). Wu et al. (2016) proposed a cyber-physical sensing and learning framework called ADDSEN to handle drone swarms for urban sensing.

We focus on the traffic monitoring application. In this problem, we assume that traffic segments that have been monitored periodically with UAVs cost less to clear incidents that subsequently occur. A recent deployment of UAVs for traffic monitoring can be seen in Figure 1, which illustrates the ongoing efforts to use UAVs in this context. The need for monitoring is assumed to be linked to traffic volumes, such that higher volumes would place a larger demand for monitoring. Realistically, however, these volumes will vary randomly over time.

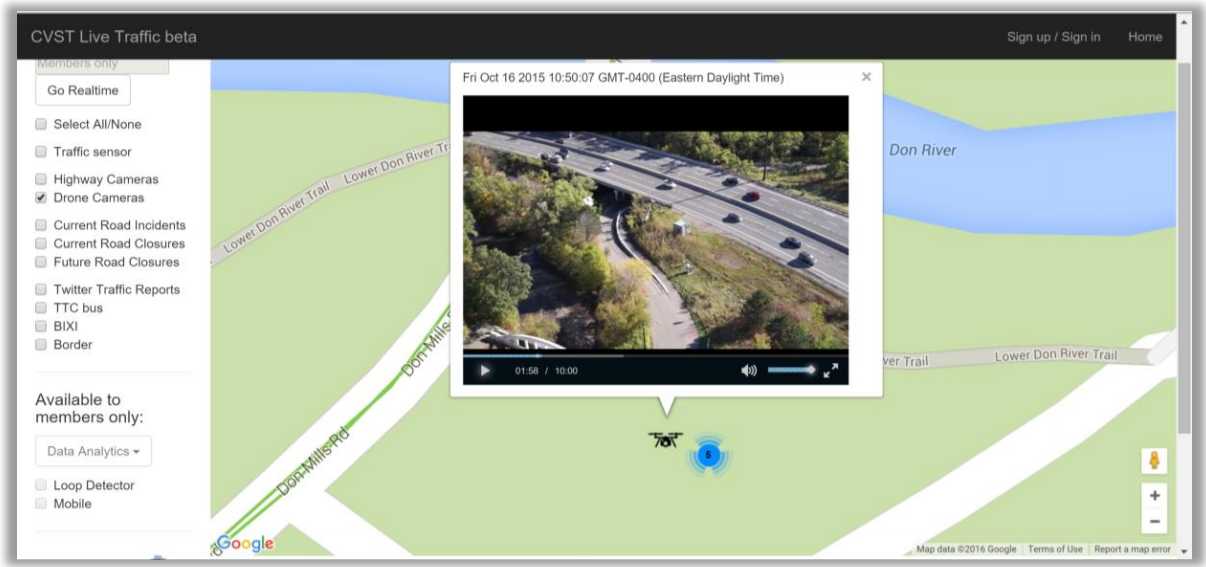

**Figure 1. Traffic data portal from Dr. Al Leon-Garcia's research group at University of Toronto (portal.cvst.ca; accessed Feb 2016).**

There are few methods of drone deployment for addressing this problem. Kinney et al. (2005) noted that current practice in planning the routes of UAVs typically involves manual calculations. More recent studies sought to address the deployment problem in several ways. Initial studies treated the monitoring problem as a multiple traveling salesman problem (mTSP) with time windows (Ryan et al., 1998; Kinney et al., 2005; Rathinam et al., 2007). Since monitoring entails repeated visits over multiple periods, researchers introduced periodic coverage and timing of mobile sensors over different areas (Cheng et al., 2011; Du et al., 2010), called "sweep coverage" (see Gage, 1992). One study looked at multi-period UAV deployment as a dynamic vehicle routing problem (Bullo et al., 2011), and another as a time-space sensor assignment problem (Zhang et al., 2015).

However, these studies all treat monitoring destinations as nodes in a network. Since UAVs are mobile sensors that can monitor while in motion, it is more appropriate to monitor arcs instead of nodes. One study focused on the arc routing aspect of mobile sensors on infrastructure networks, ensuring that tours visit each critical arc in a network like in a rural postman problem (Sipahioglu et al., 2010). Although UAVs can travel over the air and are not physically restricted to a particular network, flying over uncharted space can lead to obstructions with unknown objects. As such, arc-based traversal over an agency's own known right-of-way makes sense when considering deployment policies. In fact, non-road right of way can be set up as access links that do not need monitoring (by setting the demand for monitoring to be zero for those arcs).

Yazici et al. (2014) extended the arc routing to a dynamic deployment problem where incidents require updating the routes of the autonomous sensors in real time. While these studies considered arc routing, neither handled multiple periods.

In short, there is no one method that considers systematic UAV sensor deployment strategies with (1) arc routing or (2) periodic coverage. No study has considered UAV deployment with (3) uncertainty in demand via an online/dynamic deployment policy.

Three primary contributions are made in this study. First, we propose a new policy under uncertainty to solve this UAV traffic monitoring problem based on arc-inventory routing, where demand for monitoring is distributed over arcs instead of nodes. This is a generalized policy of which UAV deployment is one application. Other applications are also possible: mobile sensors, dynamic postal parcel deliveries, dynamic repairman problem, etc., although further customizations may be needed. As an inventory routing problem (IRP), the vehicle routing problem is operated in a multi-period setting where demand evolves as inventory levels that need restocking. The demand may be related to traffic flow densities throughout the day, for example, if real time traffic surveillance was the objective. To the best of our knowledge, there are no models of deterministic arc-inventory routing problems (AIRPs), much less for stochastic dynamic policy variations.

Second, we propose finding the policy by reformulating its Bellman equation as a "selective vehicle routing problem (VRP)" and prove its equivalence. This is important because it allows us to integrate the routing decisions directly with determining the value function.

Third, we adapt an approximate dynamic programming algorithm used in Chow and Regan (2011a) to obtain a policy. The algorithm is designed to exploit the integrated structure of the selective VRP and shown to be effective compared to myopic and static policies over several instances that are typical in size in the field of stochastic dynamic network optimization.

The remainder of this study is organized as follows. Section 2 provides a review of the literature relevant to development of the proposed policy and solution algorithm. Section 3 presents the new model by organically expanding from a core capacitated arc routing problem initially formulated by Golden and Wong (1981) to a stochastic dynamic AIRP (SDAIRP) policy designed for UAVs. A solution method is proposed for this newly formulated model and evaluated in Section 4. Section 5 presents computational experiments to evaluate the algorithm and unique qualities of the model, and Section 6 concludes.

## 2. Literature review

At the core of the problem is the need to dynamically route the mobile sensors along arcs instead of nodes. We first examine the arc routing literature.

The static routing problem for a single vehicle over arcs of a network is well known in the literature as a Chinese Postman Problem. For multiple vehicles with capacities, the problem is known as a Capacitated Arc Routing Problem (CARP), which Golden and Wong (1981) formulated as a mixed integer programming problem. They showed that this problem can be quite different from the node-based mTSP. This is because arcs may be visited more than once to reach other arcs.

As a monitoring strategy, it is also important to consider deployment over multiple periods, as suggested in some of the UAV literature (Cheng et al., 2011; Du et al., 2010; Zhang et al., 2015). For the CARP, multi-period extensions have been proposed before (Lacomme et al., 2005; Monroy et al., 2013) but they employ a frequency variable, i.e. having to visit an arc every $k$ periods. While the frequency-based approach works fine for the deterministic setting, they are not time-dependent so they don't fit into a dynamic policy framework well. For example, an accident on the road may lead to significant increases in need for monitoring over multiple upstream arcs in the network over the next two periods. A frequency based approach cannot capture this scenario. On the other hand, the scenario can be modeled using an inventory-based framework. In the example, the significant increase in need can be modeled as added "consumption" of need on each arc. Deployment of a UAV to monitor an arc would "reset" the need as an "order-up-to" inventory policy.

There are inventory routing problems (IRPs) in the literature (see Bell et al., 1983; Bertazzi et al., 2008) that integrate the inventory restocking with vehicle routing problem to minimize transport and inventory costs. However, there have been no arc-inventory routing problems that look at delivery of goods to arcs.

Furthermore, research that considers IRPs with stochastic demand remains limited. It has been studied in several cases (Bard et al., 1998; Jaillet et al., 2002; Kleywegt et al., 2004; Adelman, 2004; Hvattum et al., 2009; Bertazzi et al., 2013; Coelho et al., 2014) with approximate dynamic programming solution methods. In such methods, the policy value is defined as $V_t(\mathbb{S}_t) = \max_{a_t} \sum_{s=t}^{t+T} \gamma^s C_s(\mathbb{S}_s, a_s)$, and can modeled as a state-dependent recursive Bellman equation as illustrated in Eq. (1).

$$V_t(\mathbb{S}_t) = \max_{a_t}(C_t(\mathbb{S}_t, a_t) + \gamma E[V_{t+1}(\mathbb{S}_{t+1})|\mathbb{S}_t, a_t]) \qquad (1)$$

where $V_t(\mathbb{S}_t)$ is the value of the policy being maximized, $C_t$ is the immediate payoff at time step $t$ of the decision/action $a_t$ under state $\mathbb{S}_t$ (which is also typically driven by information on exogenous stochastic variables), and $\gamma$ is a discount factor. The challenge is in determining an appropriate value for the last term, $E[V_{t+1}(\mathbb{S}_{t+1})|\mathbb{S}_t, a_t]$, and approximation techniques are typically used.

However, these studies have all assumed time-independent stochastic variables. This unfortunately does not take full advantage of the value that real-time demand information can provide to a policy. Berman and Larson (2001) consider time-dependent stochastic demand, modeled as a mean-reverting stochastic process (see Chow and Regan, 2011b), but only apply it to a single vehicle tour.

Why should we consider using a mean-reverting process? Yang and Chu (2011) showed empirically that traffic flow data (based on the Interstate 95 highway data) can be fit with a mean-reverting process. In such a process, demand $r_{ij}(t)$ for an arc $(i,j)$ at time $t$ past an initial condition may be modeled as shown in Eq. (2).

$$dr_{ij} = \theta_{ij}(\mu_{ij} - r_{ij})dt + \sigma_{ij}dB_t \qquad (2)$$

In this equation, $\theta_{ij} > 0$ is the mean reversion rate, $\mu_{ij}$ is the mean, $\sigma_{ij} > 0$ is the volatility parameter, and $dB_t \sim N(0, t)$ is an increment in the Wiener process. Chow and Regan (2011b) showed that the process can be fitted to fire weather data to inform deployment of air tanker servers.

Can an approximate dynamic programming methodology be applied to stochastic IRP (SIRP) for multiple vehicles using real time data modeled as mean-reverting processes? Solving this problem would provide a policy for deploying UAVs in real time for city monitoring, and is also applicable to other sensor deployment strategies and stochastic inventory routing applications.

Kleywegt et al. (2004) discussed the use of value approximation functions with polynomial functions in solving their SIRP. Polynomial function-based approximation has been around since Bellman and Dreyfus (1959), who approximate functions using Legendre polynomials. More recent value function approximation methods use least squares estimation of Laguerre or Hermite polynomials to estimate the value function within a Monte Carlo simulation for handling time-dependent stochastic processes like Brownian motion or mean-reverting processes (Carriere, 1996; Longstaff and Schwartz, 2001; Gamba, 2002; Chow and Regan, 2011a, 2011c). This "Least Squares Monte Carlo" (LSM) simulation method was partially proven by Longstaff and Schwartz (2001) to converge asymptotically to the unbiased estimator of the true option value, and Stentoft (2004) further proved that the LSM approximation converges to the true value as the number of sample paths $|P| \to \infty$ if the number of polynomial basis functions $M = M(P)$ is increasing in $P$ such that $M \to \infty$ and $\frac{M^3}{|P|} \to 0$. Details of the proofs can be found in their studies, and a more comprehensive review is available from Chow and Sayarshad (2015).

None of the LSM methods have been applied to a SIRP. Chow and Regan (2011a) applied the LSM to timing a generic network design investment without any significant modifications to exploit the structure of the problem. In the case of the SIRP with mean-reverting demand processes, there is an opportunity to exploit its integration of vehicle routing with inventory control to more elegantly apply LSM.

## 3. Proposed policy
### *3.1. Problem definition*
Consider a network $G(N, A, C)$ with a set of $n$ nodes ($N$), set of undirected arcs ($A$), and a matrix of arc costs ($C$). An agency has resources (UAVs) to monitor this network over multiple periods, defined in a dynamic context as a planning horizon of periods $T = \{1, ..., |T|\}$ from the current period 0. A decision is made at the end of each period of where to deploy UAVs, including the end of current period 0.

The objective is to dynamically allocate a finite set of UAVs to links in a network that need monitoring, over multiple time periods. Links that are unmonitored over a period of time are assumed to experience a high expected incident cost $\rho$ (relative to having monitoring). The cost is measured in equivalent units of deployment costs for the UAVs—deploying more UAVs to more routes would cost more, but would monitor more links more continuously and reduce the risk of unmonitored incident costs.

The increase in need for monitoring is modeled as a stochastic process $r_{ijt}$: not all links need the same frequency of monitoring, and the need changes randomly over time as a result of such factors as traffic or weather conditions. We assume that $r_{ijt}$ is an outcome of a mean-reverting process defined in Eq (2): $dr_{ij} = \theta_{ij}(\mu_{ij} - r_{ij})dt + \sigma_{ij}dB_t$. We assume an arc has an average length of time over which unmet monitoring needs lead to the expected incident cost $\rho$; this is set as a deterministic threshold $q_{ij}$ quantified in units of accumulated need $\sum_t r_{ijt}$ to represent this average length of time. The parameters $r_{ijt}, q_{ij}$, and $\rho$ can be calibrated from observed data with and without monitoring such that $q_{ij} = 1$ if there is demand for monitoring.

The monitoring function is assumed to cost additional energy for a UAV, just like how a mobile device using video streaming would consume batteries faster than one only being used for a phone call. The additional energy is modeled by having separate traversal cost $c_{ij}$ and monitoring cost $e_{ij}$.

Under this setting, the stochastic dynamic arc inventory routing policy (SDAIRP) is a dynamic resource (UAV) allocation policy with a finite horizon, in which the optimal policy would be determined from the Bellman Equation in Eq (1): $V_t(\mathbb{S}_t) = \max_{a_t}(C_t(\mathbb{S}_t, a_t) + \gamma E[V_{t+1}(\mathbb{S}_{t+1})|\mathbb{S}_t, a_t])$. $V_t(\mathbb{S}_t)$ is the value of the policy, $\mathbb{S}_t$ is the state of the system which includes the routes of each UAV, the level of need of each arc, and current consumption rate $r_{ijt}$ for monitoring need. $C_t(\mathbb{S}_t, a_t)$ is equivalent to the allocation/transportation cost of deploying UAVs, while $a_t$ is the set of deployment routes obtained from a feasible set.

This problem can be modeled as a dynamic policy. However, since no such policies have been formulated before, we start in Section 3.2 by customizing a CARP for UAV deployment, and from there evolve the problem incrementally into the proposed policy.

### *3.2. Single period CARP*
Before we can formulate the dynamic policy, we need to first define the static problem in the context of UAV deployment. Consider first the UAV deployment problem within a single period as a CARP. We adapt the formulation from Golden and Wong (1981) as follows to fit a UAV network monitoring setting. Some variables defined earlier may be defined slightly differently here in the static context, but they will be updated in subsequent sections as we evolve the model accordingly.

Notation:
*Parameters*
$G(N, A, C)$ = network with set of $n$ nodes ($N$), set of undirected arcs ($A$), and matrix of arc costs ($C$)
$c_{ij}$ = fuel cost required to traverse arc $(i,j) \in A$
$e_{ij}$ = additional fuel cost for serving arc $(i,j) \in A$
$q_{ij}$ = demand for service at arc $(i,j) \in A$, which in the context of traffic monitoring is the cumulative need for surveillance due to increasing risk from non-coverage, and may be correlated to traffic flow, calibrated with $r_{ijt}$ so that $q_{ij} = 1$ if there is demand
$W$ = vehicle fuel capacity in units of cost

$K$ = number of available vehicles

*Decision variables*
$x_{ij}^p = 1$ if arc $(i,j) \in A$ is traversed by vehicle $p$
$l_{ij}^p = 1$ if arc $(i,j) \in A$ is serviced by vehicle $p$
$f_{ij}^p$ = flow on arc $(i,j) \in A$ by vehicle $p$

In the case of UAV deployment, path length is in terms of energy consumed by each UAV, and normally an arc is set to be needing coverage or not. Since we will expand this problem to have variable degrees of coverage, we use a binary parameter $q_{ij} \in \{0,1\}$. The mixed integer programming formulation is presented in Eq. (3).

$$\min Z_1 = \sum_{i=1}^{n} \sum_{j=1}^{n} \sum_{p=1}^{K} c_{ij} x_{ij}^p \tag{3a}$$

Subject to

$$\sum_{k=1}^{n} x_{ki}^p - \sum_{k=1}^{n} x_{ik}^p = 0, \quad i = 1, \ldots, n; p = 1, \ldots, K \tag{3b}$$

$$\sum_{p=1}^{K} (l_{ij}^p + l_{ji}^p) = q_{ij}, \quad (i,j) \in A \tag{3c}$$

$$x_{ij}^p \geq l_{ij}^p, \quad (i,j) \in A, \quad p = 1, \ldots, K \tag{3d}$$

$$\sum_{i=1}^{n} \sum_{j=1}^{n} (c_{ij} x_{ij}^p + e_{ij} l_{ij}^p) \leq W, \quad p = 1, \ldots, K \tag{3e}$$

$$\sum_{k=1}^{n} f_{ik}^p - \sum_{k=1}^{n} f_{ki}^p = \sum_{j=1}^{n} l_{ij}^p, \quad i = 2, \ldots, n; p = 1, \ldots, K \tag{3f}$$

$$f_{ij}^p \leq n^2 x_{ij}^p, \quad (i,j) \in A, \quad p = 1, \ldots, K \tag{3g}$$

$$f_{ij}^p \geq 0 \tag{3h}$$

$$x_{ij}^p, l_{ij}^p \in \{0,1\} \tag{3i}$$

where Eq. (3a) is the objective of minimizing traversal costs, Eq. (3b) conserves flow, Eq. (3c) ensures demand is met, Eq. (3d) requires an arc to be traversed before it can be covered, Eq. (3e) is the energy capacity, Eq. (3f) – (3g) cover the sub-tour elimination, and Eq. (3h) – (3i) are non-negativity and integral constraints. The model can be illustrated with a network derived from Monroy et al. (2013). We consider the 5-node graph as shown in Figure 2.

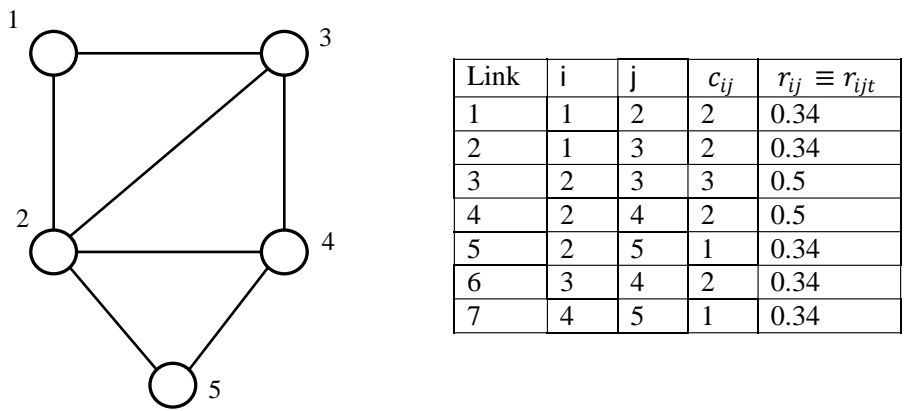

**Figure 2. Illustrative network.**

Assuming $e_{ij} = 0.1c_{ij}$, $K = 2$, and $W = 50$, the objective value of the optimal solution can be found using a commercial solver (MATLAB R2014b, *intlinprog* function) as $Z_1^* = 15$. When $W = 12$, the optimal solution changes to $Z_1^* = 19$ and requires deploying 2 UAVs. The routes are shown in Figure 3, illustrating the sensitivity of the solutions to fuel capacity parameter. For $W = 12$, some undirected arcs are visited multiple times even by the same vehicle as well as by both vehicles, simply because of the need to traverse the arc to reach another arc and limited energy constraints.

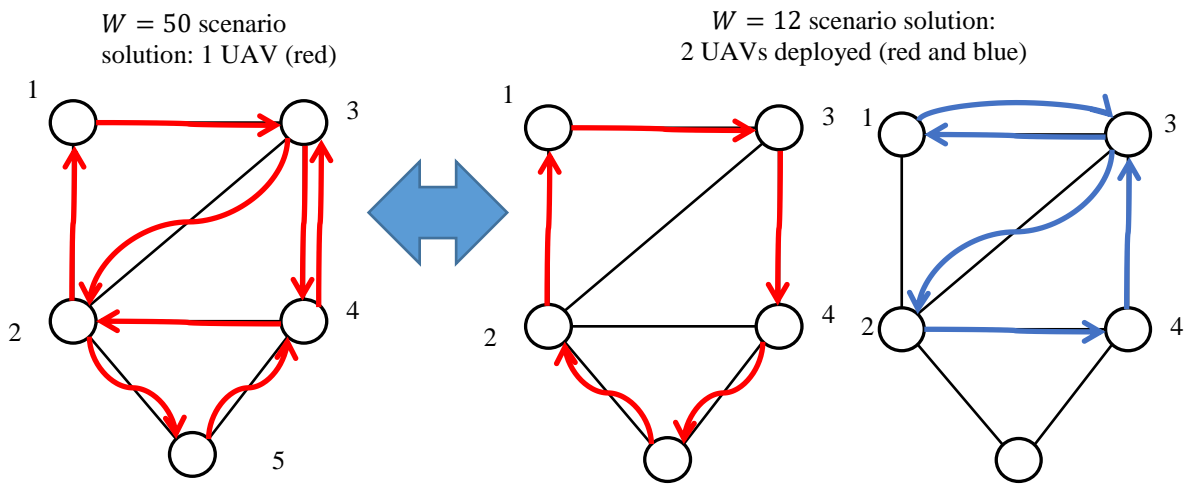

**Figure 3. Routes from $W = 50$ scenario to $W = 12$ scenario.**

### 3.3. Deterministic arc inventory routing problem for multiple periods

The first step to evolve the model in Eq (3) is to consider multiple time periods. We extend Eq (3) into an inventory routing problem for arcs, which has not been studied in the literature. The length of a time period is assumed to be long enough to allow a UAV to traverse its route completely. A new parameter $h_{ij}$ is introduced to model the inventory holding cost, which can be interpreted as cost due to overly frequent monitoring. In some cases this may be set to zero; in communities where privacy is a big issue, frequent monitoring may be a social cost. Currently, commercial UAVs can operate for a little under half an hour, so an hourly time period would make sense. Real time information like traffic data is typically collected on a more frequent basis (e.g. 5-minute intervals) so that should not be a constraining criterion. Several new variables are added to those introduced in Section 3.2.

Additional notation:
$s_{ijt}$ = unmet need for demand $q_{ij}$ at the end of time interval $t \in \{0, T\}$, where $T = \{1, ..., |T|\}$
$\zeta$ = number of periods over which a UAV needs to be recharged before it can fly again
$h_{ij}$ = inventory holding cost per period as cost of frequent monitoring
$r_{ij}$ = consumption rate of inventory per period
$\mathcal{M}$ = "big M", an arbitrarily big constant

The formulation is shown in Eq. (4).

$$\min Z_2 = \sum_{i=1}^{n}\sum_{j=1}^{n}\sum_{p=1}^{K}\sum_{t=0}^{T} c_{ij} x_{ijt}^p + \sum_{i=1}^{n}\sum_{j=1}^{n}\sum_{t=0}^{T} h_{ij} s_{ijt} \quad (4a)$$

Subject to

$$\sum_{k=1}^{n} x_{kit}^p - \sum_{k=1}^{n} x_{ikt}^p = 0, \quad i = 1,...,n, \quad p = 1,...,K, \quad t \in T \quad (4b)$$

$$\sum_{p=1}^{K} (l_{ijt}^p + l_{jit}^p) \leq 1, \quad (i,j) \in A, \quad t \in T \quad (4c)$$

$$s_{ij,t-1} - r_{ij} - s_{ijt} \geq -\mathcal{M} \sum_{p=1}^{K} (l_{ijt}^p + l_{jit}^p), \quad (i,j) \in A, \quad t \in T \quad (4d)$$

$$s_{ij,t-1} - r_{ij} - s_{ijt} \leq \mathcal{M} \sum_{p=1}^{K} (l_{ijt}^p + l_{jit}^p), \quad (i,j) \in A, \quad t \in T \quad (4e)$$

$$s_{ijt} \geq q_{ij} - \mathcal{M}\left(1 - \sum_{p=1}^{K}(l_{ijt}^p + l_{jit}^p)\right), \quad (i,j) \in A, \quad t \in T \quad (4f)$$

$$x_{ijt}^p \geq l_{ijt}^p, \quad (i,j) \in A, \quad p = 1,...,K, \quad t \in T \quad (4g)$$

$$\mathcal{M}\left(1 - \sum_{j=2}^{n} x^p_{1jt}\right) \geq \sum_{j=2}^{n} x^p_{1j,t+\tau}, \qquad p = 1, \ldots, K, \qquad t \in T, \qquad \tau = 1 \ldots \zeta, \zeta > 0 \qquad (4h)$$

$$\sum_{i=1}^{n}\sum_{j=1}^{n}\left(c_{ij}x^p_{ijt} + e_{ij}l^p_{ijt}\right) \leq W, \qquad p = 1, \ldots, K, \qquad t \in T \qquad (4i)$$

$$\sum_{k=1}^{n} f^p_{ikt} - \sum_{k=1}^{n} f^p_{kit} = \sum_{j=1}^{n} l^p_{ijt}, \qquad i = 2, \ldots, n, \qquad p = 1, \ldots, K, \qquad t \in T \qquad (4j)$$

$$f^p_{ijt} \leq n^2 x^p_{ijt}, \qquad (i,j) \in A, \qquad t \in T, \qquad p = 1, \ldots, K \qquad (4k)$$

$$f^p_{ijt} \geq 0 \qquad (4l)$$

$$r_{ij} \leq s_{ijt} \leq q_{ij}, \qquad (i,j) \in A, \qquad t \in \{0, T\} \qquad (4m)$$

$$x^p_{ijt}, l^p_{ijt} \in \{0,1\} \qquad (4n)$$

Objective (4a) is divided into two minimizing terms: the energy consumption representing transportation cost, and an inventory holding cost. Constraints (4b), (4c), (4g), (4i), (4j), (4k), (4l) and (4n) are mostly the same as Eq. (3), except the variables are expanded into multiple periods. Constraints (4d) – (4f) update the fuel from one period to the next to either decrease if there is no monitoring, or to increase up to the maximum if monitored in a period. This constraint is similar to the refueling constraint in Chow and Liu (2012). Constraint (4h) ensures that a UAV that is refueling is out of commission for $\zeta$ periods. Constraints (4m) are the inventory bound constraints to prevent the inventory from dropping below zero or going above the maximum $q_{ij}$. Although this problem introduces a time dimension, it is a multi-period static problem for the UAVs (while demand is dynamic), in the sense that the trip lengths of the UAVs are not tracked over multiple periods. For dynamic UAV tracking, Zhang et al. (2015) provide a clever time-geographic approach to model it.

We test this formulation out with a commercial solver (*intlinprog* on MATLAB) on the same example shown in Figure 2, with $W = 12$ and $h_{ij} = 0.1 \, \forall (i,j) \in A$. The optimal solution is $Z_2 = 41.07$ (run time 5866.36 sec on 64-bit Windows 7 SP1 Intel Core i7-3770 CPU with 3.40GHz and 16GB RAM), and the solution routes are shown in Figure 4.

### 3.4. Stochastic dynamic arc inventory routing policy (SDAIRP)

As mentioned in Section 3.1, we assume that the incremental need for monitoring is directly related to random traffic conditions, i.e. $r_{ijt}$ is assumed to be an outcome of a mean-reverting process defined in Eq. (2): $dr_{ij} = \theta_{ij}(\mu_{ij} - r_{ij})dt + \sigma_{ij}dB_t$. Arc monitoring decisions ($y_t$) are made at time $t$ with current information on the need for monitoring ($s_t$) and historical information on rate of change ($r_t$). The monitoring decision is assumed to occur after updating the change, so a stock out can occur if $s_t - r_t < 0$. Different policies can be applied, as illustrated in Figure 5 for a single stochastic arc.

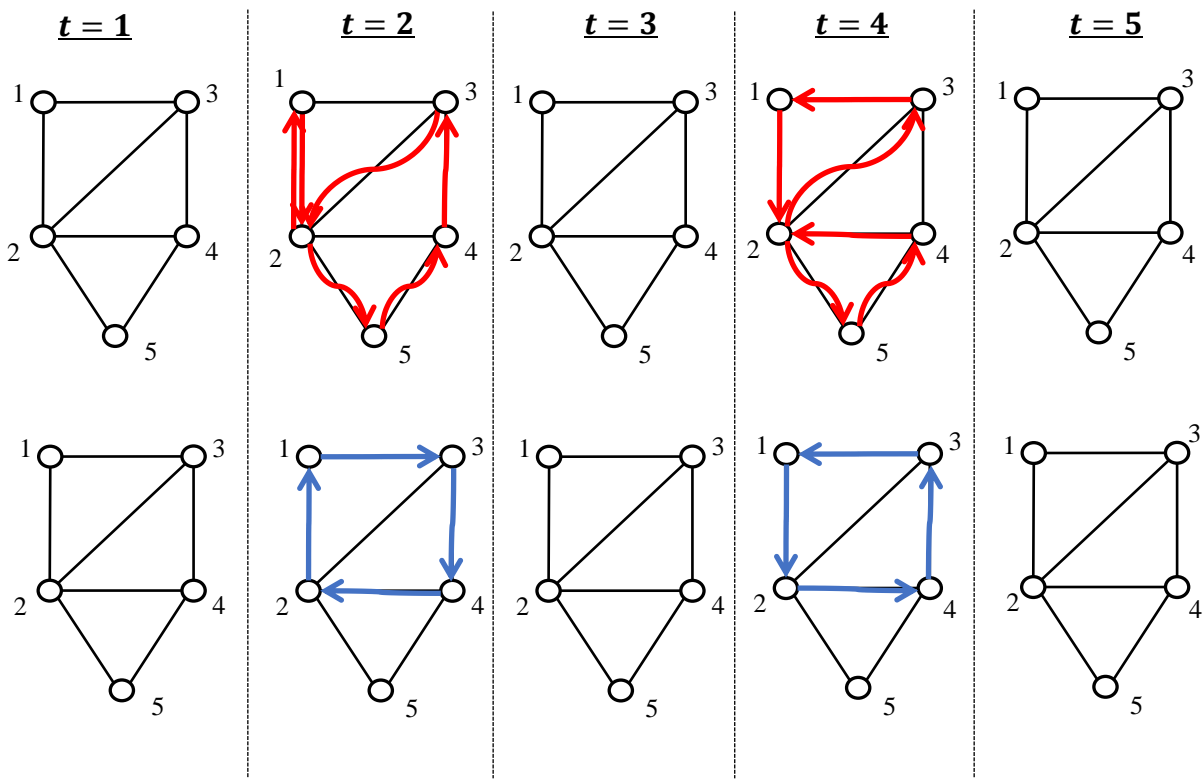

**Figure 4.** Routes from $W = 12$ scenario for arc-inventory routing problem over 5 periods.

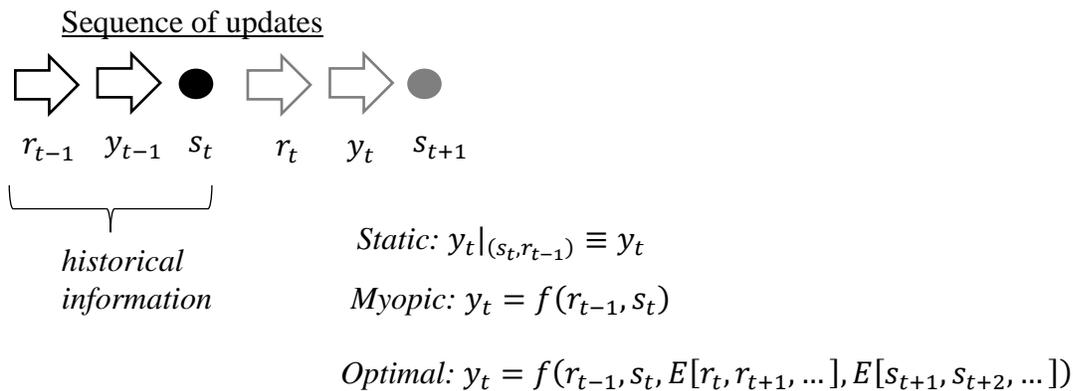

**Figure 5.** Illustration of different dynamic policies in this setting.

The SDAIRP has an optimal policy that can be determined dynamically from the Bellman Equation in Eq (1): $V_t(\mathbb{S}_t) = \max_{a_t}(C_t(\mathbb{S}_t, a_t) + \gamma E[V_{t+1}(\mathbb{S}_{t+1})|\mathbb{S}_t, a_t])$. $V_t(\mathbb{S}_t)$ is the value of the policy, $\mathbb{S}_t$ is the state of the system which includes the routes of each UAV, the level of need of each arc, and current consumption rate $r_{ijt}$ for monitoring need. $C_t(\mathbb{S}_t, a_t)$ is equivalent to the allocation/transportation cost of deploying UAVs, while $a_t$ is the set of deployment routes obtained from a feasible set.

The challenge of solving Eq. (1) is the curse of dimensionality in the problem, since $r_{ij}$ is a mean-reverting process defined over a continuous domain rather than several discrete state values. We seek a simplification to obtain the policy value $V_t(\mathbb{S}_t)$ by approximating the future time horizon with only the time until the next replenishment $t+\tau$ with $\hat{V}_t(\mathbb{S}_t) = \max_{a_t}\sum_{s=t}^{t+\tau}\gamma^s C_s(\mathbb{S}_s, a_s)$, as shown in Figure 6. This approximation is justified by the order-up-to policy which resets the inventory level at replenishment, reducing the error that may occur from such an approximation.

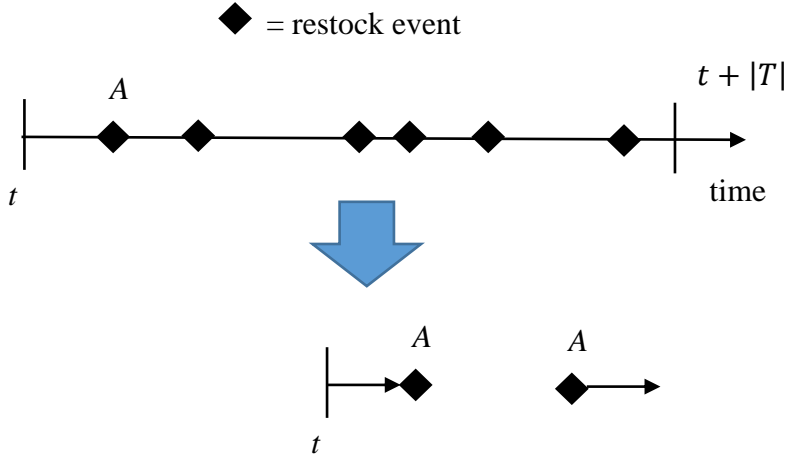

**Figure 6. Illustration of inventory replenishment with random consumption under a stationary stochastic process.**

We can exploit the structure of the inventory routing problem to obtain a policy that is optimal to that approximation. This is presented in Proposition 1.

**Proposition 1.** *Consider the stochastic dynamic arc inventory routing policy with order-up-to replenishment. For its approximation $\hat{V}_t(\mathbb{S}_t) = \max_{a_t}\sum_{s=t}^{t+\tau}\gamma^s C_s(\mathbb{S}_s, a_s)$ until replenishment time $t + \tau$, a selective vehicle routing problem can be used recursively to find its optimal value.*

**Proof**. Consider the series of decisions being made in each time increment (up until the next replenishment, where $(i,j)$ notation is dropped for simplicity):

$$\hat{V}_t(\mathbb{S}_t) = \max_{a_t}(C_t(\mathbb{S}_t, a_t) + \gamma E_\tau[V_{t+1}(\mathbb{S}_{t+1})|\mathbb{S}_t, a_t])$$

where $E_\tau$ is the expected value up until the next replenishment occurring at time $t + \tau$. The approximation is a stopping problem, where the latter term $E_\tau[V_{t+1}(\mathbb{S}_{t+1})|\mathbb{S}_t, a_t] = 0$ if the site is covered and inventory is replenished. If the site is not covered, the inventory would decrease by a random amount $r_t$ until the next time increment $t+1$. Since there is no routing occurring at a time $t$ if coverage is deferred, the value $V_t(\mathbb{S}_t)$ is essentially comprised of the latter term at the time of replenishment, $E_\tau[V_\tau(\mathbb{S}_\tau)|\mathbb{S}_t]$. In other words, the payoffs are structured as:

$$y_t = 1, C_t(\mathbb{S}_t, a_t) > 0, E_\tau[V_{t+1}(\mathbb{S}_{t+1})|\mathbb{S}_t, a_t] = 0$$

$$y_t = 0, C_t(\mathbb{S}_t, a_t) = 0, E_\tau[V_{t+1}(\mathbb{S}_{t+1})|\mathbb{S}_t, a_t] > 0$$

This means the $\hat{V}_t(\mathbb{S}_t)$ can be expressed as an optimization problem that allocates fully between $C_t(\mathbb{S}_t, a_t)$ (cover now) and $\gamma E_\tau[V_{t+1}(\mathbb{S}_{t+1})|\mathbb{S}_t, a_t]$ (cover later), subject to feasibility constraints. This is in fact a selective vehicle routing problem (e.g. Gribhovskaia et al., 2008; Valle et al., 2009; Chow, 2014; Allahviranloo et al., 2014) which selects a subset of elements in a network to visit for some payoff. ∎

Let us now formulate this selective VRP.

Additional notation from Section 3.3:
$y_{ijt}$ is 1 if arc $(i,j)$ is selected to be monitored at time $t$, 0 otherwise
$\rho$ is cost of a stock out representing expected incident cost with no monitoring
$R(s_{ijt})$ is the probability of an eventual stock out given inventory $s_{ijt}$, where $R(s_{ijt})\rho = E_\tau[V_\tau(\mathbb{S}_\tau)|\mathbb{S}_t]$

The following optimization problem in Eq. (5) (the period $t$ in the decision variables is removed for simplicity in presentation) is solved recursively.

$$\min_{x,l,f,y} Z_3 = \sum_{i=1}^{n}\sum_{j=1}^{n}\sum_{p=1}^{K} c_{ij}x_{ij}^p + \sum_{i=1}^{n}\sum_{j=1}^{n}\left(\rho_{ij}R(q_{ij}) - \rho_{ij}R(s_{ijt}) + h_{ij}(q_{ij} - s_{ijt})\right)y_{ij} \quad (5a)$$

Subject to

$$\sum_{k=1}^{n} x_{ki}^p - \sum_{k=1}^{n} x_{ik}^p = 0, \quad i = 1,\ldots,n; p = 1,\ldots,K \quad (5b)$$

$$\sum_{p=1}^{K}(l_{ij}^p + l_{ji}^p) = y_{ij}, \quad (i,j) \in A \quad (5c)$$

$$x_{ij}^p \geq l_{ij}^p, \quad (i,j) \in A, \quad p = 1,\ldots,K \quad (5d)$$

$$\sum_{i=1}^{n}\sum_{j=1}^{n}(c_{ij}x_{ij}^p + e_{ij}l_{ij}^p) \leq W, \quad p = 1,\ldots,K \quad (5e)$$

$$\sum_{k=1}^{n} f_{ik}^p - \sum_{k=1}^{n} f_{ki}^p = \sum_{j=1}^{n} l_{ij}^p, \quad i = 2,\ldots,n; p = 1,\ldots,K \quad (5f)$$

$$f_{ij}^p \leq n^2 x_{ij}^p, \quad (i,j) \in A, \quad p = 1,\ldots,K \quad (5g)$$

$$f_{ij}^p \geq 0 \tag{5h}$$

$$x_{ij}^p, l_{ij}^p, y_{ij} \in \{0,1\} \tag{5i}$$

Eq. (5) is a selective vehicle routing problem where an operator decides which destination arcs to visit such that the profit is the change in expected cost of stock out from $\rho_{ij}R(s_{ijt})$ to $\rho_{ij}R(q_{ij})$, plus changes in inventory cost $h_{ij}(q_{ij} - s_{ijt})$. Eq. (5a) uses a binary selection variable $y_{ij}$ to determine if an arc $(i,j) \in A$ is covered in that period. The consequence of not selecting an arc is the increased risk of stock out compared to the case of having full stock available. Eq. (5c) requires the selection variable to be one if covered. The other constraints are mostly the same as in Eq. (3).

Unlike the deterministic arc inventory routing problem (AIRP) formulated over multiple periods, the problem in Eq. (5) is solved within each period so its computational complexity is almost the same as the single period CARP (additional binary selection variables). In terms of run time, Eq. (5) can be solved for the example in Figure 1 in 0.12 s under the computational setting.

Due to the presence of the stochastic element $E_\tau[V_\tau(\mathbb{S}_\tau)|\mathbb{S}_t]$ embedded in the $R(s_{ijt})$, running the problem recursively would still have the dimensionality problems as Eq. (1). However, an effective solution algorithm from the real option literature can be used with Eq. (5) as a sub-problem, which is presented and evaluated in Section 4.

## 4. Proposed solution algorithm
### *4.1. Least squares approximation of stock out cost*
A solution algorithm is needed to approximate the future value $R(s_{ijt})$ and derive the policy. We propose a variation of the LSM algorithm applied to the approximate SDAIRP in a finite horizon with $\gamma = 1$ and $r_{ijt}$ distributed as a mean-reverting process as shown in Eq. (2). Since the algorithm is a numerical method, the policy is not restricted to only mean-reverting processes; any other stationary stochastic process can also be applied.

In the LSM, the value function one time step ahead ($E[V_{t+1}]$) is approximated using least squares estimation of a polynomial basis function fitted to $P$ Monte Carlo-simulated independent sample paths. Longstaff and Schwartz (2001) discuss using either of two different polynomial functions, Laguerre or Hermite polynomials. Chow and Regan (2011a) employ Hermite polynomials, which can be constructed as shown in Eq. (6).

$$H_n(x) = (-1)^n e^{\frac{x^2}{2}} \frac{d^n}{dx^n} e^{-\frac{x^2}{2}} \tag{6}$$

For example, a Hermite polynomial function with $M = 5$ basis functions would have a function:

$$f(x) = \beta_0 + \beta_1 x + \beta_2(x^2 - 1) + \beta_3(x^3 - 3x) + \beta_4(x^4 - 6x^2 + 3) + \beta_5(x^5 - 10x^3 + 15x)$$

Instead of estimating one time step ahead, we modify the least squares approximation to use the **sampled stock out cost at event of replenishment** as the dependent variable. This is a fundamental change in the LSM algorithm (instead of using the sampled value at the next time step), but one that makes intuitive sense based on Proposition 1, and verifiable with a simple example.

Consider a simple myopic replenishment policy for a single location with mean-reverting process and inventory parameters $\mu = 0.5, \theta = 0.1, \sigma = 0.1, r_0 = 0.33, s_0 = 1$. If we simulate the consumption rate over 4 time steps across 12 independent sample paths, and assume restocking occurs whenever the inventory drops below 0.33, we get the following sample Table 1.

In this table, the R{t} is the consumption at time $t$, the S{t}PRE stands for the pre-replenishment decision inventory level at time $t$, while S{t}POST stands for the post-replenishment decision. The S.O. is the stock out cost at time step 4, where any negative values (the last two samples) result in paying a stock out cost of $\rho = 10$.

Table 1. Example of simulated sample paths for LSM with stock out cost as dependent variable

|    | R1    | S1    | R2    | S2PRE | S2POST | R3    | S3PRE | S3POST | R4    | S4PRE  | S.O. |
|----|-------|-------|-------|-------|--------|-------|-------|--------|-------|--------|------|
| 1  | 0.387 | 0.613 | 0.331 | 0.282 | 1.000  | 0.276 | 0.724 | 0.724  | 0.269 | 0.454  | 0    |
| 2  | 0.235 | 0.765 | 0.437 | 0.328 | 1.000  | 0.412 | 0.588 | 0.588  | 0.432 | 0.155  | 0    |
| 3  | 0.288 | 0.712 | 0.288 | 0.424 | 0.424  | 0.311 | 0.113 | 1.000  | 0.360 | 0.640  | 0    |
| 4  | 0.232 | 0.768 | 0.180 | 0.588 | 0.588  | 0.283 | 0.305 | 1.000  | 0.263 | 0.737  | 0    |
| 5  | 0.295 | 0.705 | 0.259 | 0.445 | 0.445  | 0.321 | 0.124 | 1.000  | 0.296 | 0.704  | 0    |
| 6  | 0.453 | 0.547 | 0.503 | 0.043 | 1.000  | 0.410 | 0.590 | 0.590  | 0.464 | 0.126  | 0    |
| 7  | 0.180 | 0.820 | 0.081 | 0.739 | 0.739  | 0.162 | 0.578 | 0.578  | 0.255 | 0.323  | 0    |
| 8  | 0.207 | 0.793 | 0.107 | 0.686 | 0.686  | 0.139 | 0.547 | 0.547  | 0.141 | 0.406  | 0    |
| 9  | 0.374 | 0.626 | 0.436 | 0.191 | 1.000  | 0.386 | 0.614 | 0.614  | 0.295 | 0.318  | 0    |
| 10 | 0.340 | 0.660 | 0.375 | 0.285 | 1.000  | 0.454 | 0.546 | 0.546  | 0.479 | 0.067  | 0    |
| 11 | 0.463 | 0.537 | 0.394 | 0.142 | 1.000  | 0.560 | 0.440 | 0.440  | 0.522 | -0.082 | 10   |
| 12 | 0.426 | 0.574 | 0.537 | 0.037 | 1.000  | 0.590 | 0.410 | 0.410  | 0.672 | -0.263 | 10   |

The least squares estimation is illustrated for time step 3 over the 12 samples. Using S3POST as the independent variable "x" in Eq. (6) with $M = 5$, the following Table 2 can be constructed for the Hermite polynomials.

Table 2. Hermite polynomial regression at time step 3 with $M = 5$

| M0 | M1 | M2 | M3 | M4 | M5 | OBSERVED S.O. | PREDICTED S.O. |
|---|---|---|---|---|---|---|---|
| 1 | 0.724 | -0.476 | -1.792 | 0.132 | 7.264 | 0 | -0.019 |
| 1 | 0.588 | -0.655 | -1.560 | 1.047 | 6.856 | 0 | -0.253 |
| 1 | 1.000 | 0.000 | -2.000 | -2.000 | 6.000 | 0 | 0.000 |
| 1 | 1.000 | 0.000 | -2.000 | -2.000 | 6.000 | 0 | 0.000 |
| 1 | 1.000 | 0.000 | -2.000 | -2.000 | 6.000 | 0 | 0.000 |
| 1 | 0.590 | -0.652 | -1.565 | 1.031 | 6.869 | 0 | -0.209 |
| 1 | 0.578 | -0.666 | -1.541 | 1.108 | 6.802 | 0 | -0.350 |
| 1 | 0.547 | -0.701 | -1.477 | 1.295 | 6.617 | 0 | 0.198 |
| 1 | 0.614 | -0.623 | -1.610 | 0.882 | 6.981 | 0 | 0.441 |
| 1 | 0.546 | -0.702 | -1.475 | 1.301 | 6.609 | 0 | 0.245 |
| 1 | 0.440 | -0.806 | -1.236 | 1.874 | 5.768 | 10 | 9.900 |
| 1 | 0.410 | -0.832 | -1.160 | 2.021 | 5.469 | 10 | 10.047 |
| $\beta_0$ | $\beta_1$ | $\beta_2$ | $\beta_3$ | $\beta_4$ | $\beta_5$ | | |
| -604467 | 1373749 | -1076649 | 690359 | -158836 | 48961 | | |

We run the least squares estimation with several other parameters, for $M = 2$ and $M = 4$, and plot the three sets of predicted points in Figure 6.

There are two important points to note from this replicable example. First, the Hermite polynomials do appear to be more accurate as the number of basis functions go up (based on the posted $R^2$ values in Figure 6 and from visual observation). Second, the premium due to volatility in the consumption rate can be observed. Because $\mu = 0.5$, the stock out cost should typically occur when the current inventory at time step 3 is below $s_3 = 0.5$. However, the volatility leads to a probabilistic stock out cost even for values slightly higher than 0.5, as shown in Figure 7. This result demonstrates that the method is working as intended.

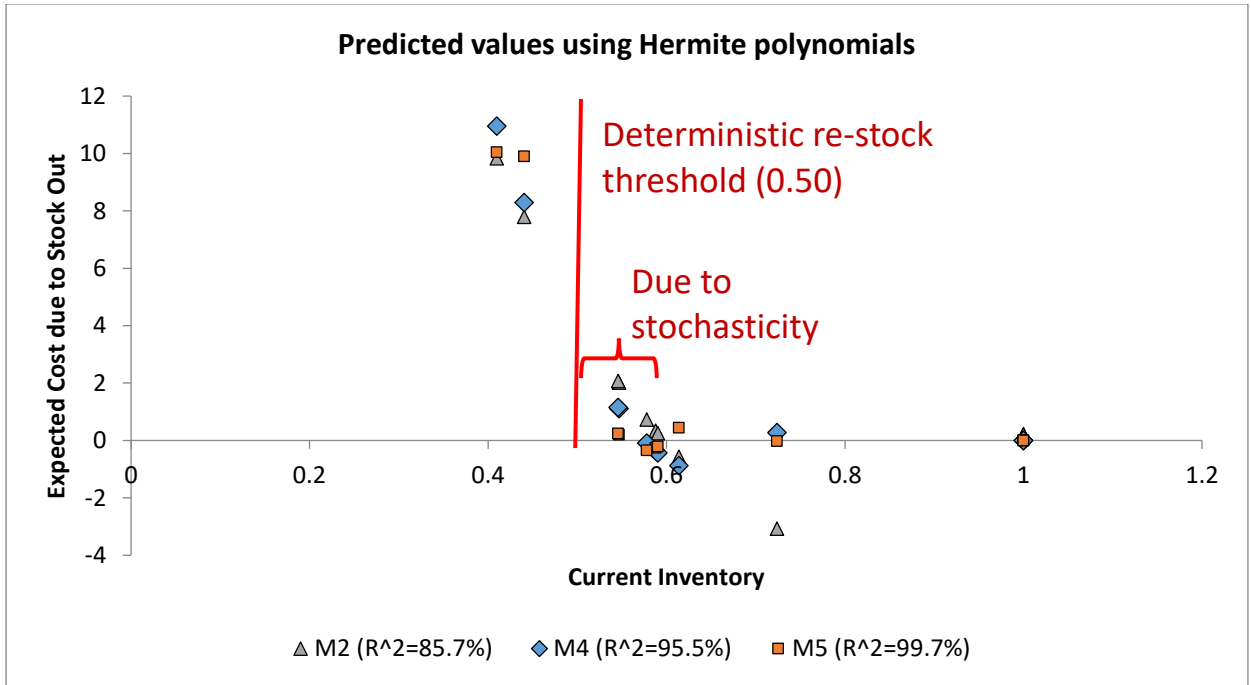

**Figure 7. Estimated stock out costs using Hermite polynomials with $M = 2, 4, 5$.**

### *4.2. Proposed algorithm*
We now introduce the proposed Algorithm 1 based on the LSM method used in Chow and Regan (2011a), but modified to regress with stock out cost dependent variables and Eq. (5) as the recursive value function.

### **Algorithm 1**
Given: initial observed consumption rate $r_{ijt}$, initial observed inventory levels $s_{ijt}$, mean-reverting process parameters $(\mu_{ij}, \theta_{ij}, \sigma_{ij})$, finite time horizon $T$, stock out cost $\rho_{ij}$, fleet size $K$, number of sample paths $P$, number of Hermite polynomial basis functions $M$, monitoring cost $e_{ij}$, undirected graph $G(N, A, C)$

1. Simulate $P$ sample paths.

a. For each simulated O-U path (see Chow & Regan, 2011b): $r_{ij,t+\Delta t} = r_{ij,t}e^{-\theta_{ij}\Delta t} + \mu_{ij}(1 - e^{-\theta_{ij}\Delta t}) + \sigma_{ij}N_{0,1}\sqrt{\frac{1-e^{-2\theta_{ij}\Delta t}}{2\theta_{ij}}}$, where $N_{0,1}$ is a random draw of standard normal distribution.
2. Solve *a priori* IRP policy $\Phi_0$ based on deterministic arc inventory routing policy or assume a naïve policy (some cyclic or periodic assignment of UAVs).
3. Update inventories $s^0_{ijp,t+\tau}$ based on $\Phi_0$ for each sample path and time step $\tau \leq |T|$.
4. For $\tau = |T|: 1$
   a. For $p = 1: P$
      i. Estimate $o_{ij}R(q_{ij}) - o_{ij}R(s_{ijt})$ using Hermite polynomial least squares based on stock out cost at "next replenishment" event. For $\tau = |T|$ this value is set to zero.
      ii. Solve a selective VRP (Eq. (5)) with the approximated costs.
      iii. Update
         1. Post-decision inventories of each link at each $(p, t + \tau)$; if replenished, then update the "next replenishment" event time
         2. Stock out occurrence and cost
5. Update policy $\Phi$ at time $t$

Output: policy $\Phi$, decisions $(x_{ijt}, y_{ijt}) = \Phi(r_{ijt}, s_{ijt})$

As the structure of the algorithm is essentially the same as the original LSM in Longstaff and Schwartz (2001), the algorithm should converge to the true value as the number of sample paths $|P| \to \infty$, where the number of polynomial basis functions $M = M(P)$ is increasing in $P$ such that $M \to \infty$ and $\frac{M^3}{|P|} \to 0$ (Stentoft, 2004).

Since the algorithm simulates independent sample paths to use for least squares, the computational cost is linearly proportional to the complexity of the underlying selective VRP in Eq. (5), which is known to be NP-hard (see Allahviranloo et al., 2014). The focus of our contribution is the dynamic policy itself, which is computationally tractable.

## 5. Computational experiments

Two sets of experiments are conducted. The first (Sections 5.1 and 5.2) is based on the same Monroy network in Figure 1 to evaluate different parameter settings for the SDAIRP. The second (Section 5.3) is based on benchmark instances from the literature to examine the scalability of the algorithm and to illustrate the dynamic process under a larger instance. Since the performance of Algorithm 1 is to be evaluated, all the evaluations use exact algorithms for the VRP sub-problems (Step 4.ii in Algorithm 1) to avoid adding unnecessary noise.

### *5.1. Experimental design*
In this experiment, we test the algorithm on the same network in Figure 1, applied sequentially over 5 time periods and 500 sample paths (resulting in 2500 runs of selective VRP for a single initial condition, resulting in 302,500 runs over 30 simulated trajectories with 4 decision points in each plus one initial decision). The experimental design is as follows.

The purpose of this experiment is to determine whether the proposed policy under different parametric designs can outperform two other benchmark policies: a static policy and a myopic policy that replenishes every time the inventory drops below a set threshold ($\mu_{23} = 0.5$). As this is a new problem, there are no other solution algorithms to compare the proposed algorithm against without significant new algorithm development. Upon publication, the data and code will be uploaded to an online repository to serve as benchmark performance measures for other algorithms (such as policy iteration and value iteration methods) to be compared. Furthermore, the results against a myopic policy can in turn be used for other algorithms developed and compared to myopic policies, with some caveats (see Chow and Sayarshad, 2015).

Policies:
- Static policy – using the deterministic AIRP solution in Figure 4
- Myopic policy – replenishing if the inventory level is $s_{ijt} < 0.5$
- Proposed SDAIRP policy, with three variations:
  - $T = 2, M = 5$
  - $T = 5, M = 5$
  - $T = 5, M = 10$

To clearly distinguish the factors to the performances, only one arc is treated as a stochastic variable, link (2,3), with parameters set at time $t = 0$ to be $\mu_{23} = 0.5, \theta_{23} = 0.1, \sigma_{23} = 0.1, r_{23,0} = 0.5$. The other arcs evolve deterministically. Additional parameters are shown in Table 3.

**Table 3. Additional parameters for the stochastic dynamic policy evaluation**

| Link | i | j | $s_{ij,0}$ |
|---|---|---|---|
| 1 | 1 | 2 | 0.68 |
| 2 | 1 | 3 | 0.68 |
| 3 | 2 | 3 | 1 |
| 4 | 2 | 4 | 1 |
| 5 | 2 | 5 | 0.68 |
| 6 | 3 | 4 | 0.68 |
| 7 | 4 | 5 | 0.68 |

Performance is measured by the simulating the $r_{23}$ over 5 time periods for 30 realizations. For example, for one realization at time period 3, the static policy would use the solution given for time period 3, the myopic policy would check the pre-decision $s_{23,3}$, and the proposed policy (with $T = 5$) would simulate 500 sample paths from $t = 3$ to $t = 8$ to determine the policy at $t = 3$. With the policy decisions simulated for the 30 trajectories, average values of total travel cost (X), total inventory holding cost (H), and total stock out cost (O) can be aggregated for each period and in total.

### 5.2. Monroy network results

The averaged results over the 30 realizations for the five test policies are shown in Table 4. Findings from this experiment are highlighted in bullet form.
- The worst performing policy is the myopic policy in this example, even compared to the static policy. This type of result has been shown in Chow and Regan (2011b), where high transportation cost and low stock out cost, along with a short time frame (only 5 periods) and high consumption rate (average of 2-3 periods to stock out) can result in the myopic policy over spending on the safe side. This is verified by the lack of stock outs, but the heightened cost of transport.

- The proposed policy outperforms both myopic and static policies in all three parametric designs.
- The two-period horizon ($T = 2$) policy outperforms the five-period horizon policy in this example. This is likely due to the high consumption rate resulting in stock out or replenishment every 2-3 periods.
- The policy with higher number of basis functions ($M = 10$) outperforms the base proposed policy, which confirms that increasing the number of basis functions improves the estimation of the value function.
- It is notable that the 23% to 28% performance improvements over myopic policy are due only to fluctuations in *one* out of seven links. This suggests that additional uncertainty in the system could further magnify the performance improvement.

**Table 4. Comparison of average (over 30 realizations) costs accumulated via different test policies**

| Static | X | H | O | Total (vs Myopic) |
|---|---|---|---|---|
| t=1 | 0.0000 | 0.2703 | 0.0000 | 0.2703 |
| t=2 | 19.000 | 0.7000 | 0.0000 | 19.7000 |
| t=3 | 0.0000 | 0.4268 | 0.0000 | 0.4268 |
| t=4 | 19.000 | 0.7000 | 0.0000 | 19.7000 |
| t=5 | 0.0000 | 0.4245 | 0.0000 | 0.4245 |
| Total | 38.0000 | 2.5215 | 0.0000 | 40.5215 (-17.5%) |
| **Myopic** | X | H | O | Total |
| t=1 | 0.0000 | 0.2703 | 0.0000 | 0.2703 |
| t=2 | 19.000 | 0.7000 | 0.0000 | 19.7000 |
| t=3 | 4.2000 | 0.4658 | 0.0000 | 4.6658 |
| t=4 | 18.4000 | 0.6910 | 0.0000 | 19.0910 |
| t=5 | 4.9000 | 0.4682 | 0.0000 | 5.3682 |
| Total | 46.5000 | 2.5953 | 0.0000 | 49.0953 |
| **SDAIRP (T=2, M=5)** | X | H | O | Total |
| t=1 | 0.0000 | 0.2703 | 0.0000 | 0.2703 |
| t=2 | 3.9667 | 0.0642 | 0.3333 | 4.3642 |
| t=3 | 17.3000 | 0.6651 | 0.0000 | 17.9651 |
| t=4 | 3.2667 | 0.4540 | 0.0000 | 3.7207 |
| t=5 | 8.4667 | 0.4030 | 0.3333 | 9.2030 |
| Total | 33.000 | 1.8566 | 0.6667 | 35.5232 (-27.6%) |
| **SDAIRP (T=5, M=5)** | X | H | O | Total |
| t=1 | 0.0000 | 0.2703 | 0.0000 | 0.2703 |
| t=2 | 7.3000 | 0.1940 | 0.3333 | 7.8274 |
| t=3 | 13.7000 | 0.6168 | 0.0000 | 14.3168 |
| t=4 | 5.2667 | 0.4577 | 0.0000 | 5.7244 |
| t=5 | 8.8000 | 0.4401 | 0.0000 | 9.2401 |
| Total | 35.0667 | 1.9790 | 0.3333 | 37.3790 (-23.9%) |
| **SDAIRP (T=5, M=10)** | X | H | O | Total |
| t=1 | 0.0000 | 0.2703 | 0.0000 | 0.2703 |
| t=2 | 7.8667 | 0.2310 | 0.3333 | 8.4310 |
| t=3 | 12.7333 | 0.6035 | 0.0000 | 13.3368 |
| t=4 | 5.3333 | 0.4553 | 0.0000 | 5.7886 |
| t=5 | 7.2000 | 0.4040 | 0.3333 | 7.9373 |
| Total | 33.1333 | 1.9641 | 0.6667 | 35.7641 (-27.2%) |

These results taken over 30 different sample instances and five time periods give empirical evidence for the effectiveness of the proposed policy.

*5.3. Benchmark instances*

In the second set of experiments, two objectives are set. The first one deals with solving the arc routing problem on progressively larger instances to provide benchmarks under the exact solution approach using commercial software. The second objective is to run the algorithm for a larger instance where all links are stochastic and to illustrate the dynamic updating of the algorithm.

The classic benchmark test instances used are the *gdb19*, *gdb15*, and *gdb9* from Golden et al. (1983). The characteristics of these instances, and their computational performances, are shown in Table 5 along with the Monroy network results from Figure 1. In all instances a feasible solution was found. It is clear from this test that network instances with more than 20 nodes using exact algorithms from commercial software would not be recommended.

**Table 5. Running CARP exact algorithms from commercial software for benchmark instances**

| Instance | No. nodes | No. links | No. vehicles | Vehicle capacity | Computed Results | |
| --- | --- | --- | --- | --- | --- | --- |
| | | | | | Run time for MATLAB *intlinprog* (sec) | No. nodes explored in algorithm |
| **Monroy (for comparison)** | 5 | 7 | 2 | 50 | 0.10 | 635 |
| *gdb19* | 8 | 11 | 3 | 27 | 1.08 | 5306 |
| *gdb15* | 7 | 21 | 4 | 37 | 4847.40 | 10000000 (max) |
| *gdb9* | 27 | 51 | 10 | 27 | 7200+ (timed out) | -- |

Practitioners can either: (a) run heuristics in the VRP sub-problem to obtain reasonable but noisy solutions; (b) break up a city-scale region into districts that each encompass 10-15 significant corridors to monitor; or (c) consider continuous approximation schemes for the underlying VRP. Districting has been shown in the literature (e.g. Stein, 1978) to be an effective practical solution for tackling large scale routing problems.

To further illustrate the use of the proposed algorithm as a dynamic, sequential policy, the *gdb19* instance is expanded upon. Figure 8a illustrates a not-to-scale visualization of the network. Figure 8b lists the additional attributes added to *gdb19* to create a new instance *gdb19-sdairp*.

The proposed algorithm is run on a simulated trajectory of $r_{ij}$'s over 15 periods to illustrate how it works when all links are stochastic. For this instance, a stock out cost of 100 units is assumed, and a finite horizon of $T = 5$ periods, $P = 100$ sample paths, and $M = 10$ basis functions are assumed for the algorithm. The realized changes in need for monitoring due to the algorithm are shown in Figure 9.

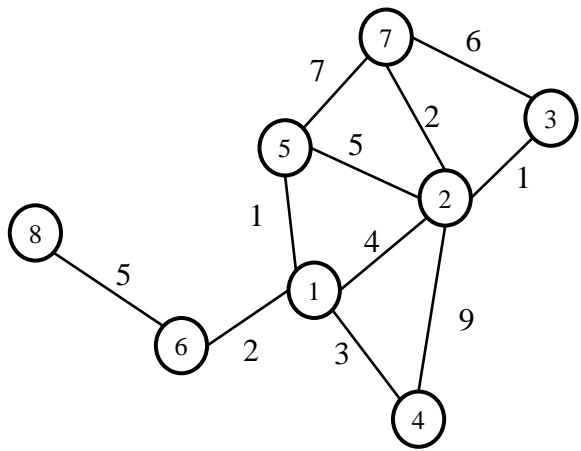

(a)

| Link | $\mu_{ij}$ | $\theta_{ij}$ | $\sigma_{ij}$ | $s_{ij,0}$ | $r_{ij,0}$ |
|---|---|---|---|---|---|
| (1,2) | 0.40 | 0.20 | 0.039 | 0.99 | 0.40 |
| (1,4) | 0.15 | 0.29 | 0.003 | 0.96 | 0.15 |
| (1,5) | 0.25 | 0.34 | 0.006 | 0.57 | 0.25 |
| (1,6) | 0.40 | 0.21 | 0.018 | 0.87 | 0.40 |
| (2,3) | 0.20 | 0.15 | 0.013 | 0.62 | 0.20 |
| (2,4) | 0.30 | 0.29 | 0.002 | 0.34 | 0.30 |
| (2,5) | 0.05 | 0.33 | 0.004 | 0.46 | 0.05 |
| (2,7) | 0.45 | 0.17 | 0.033 | 0.84 | 0.45 |
| (5,7) | 0.45 | 0.11 | 0.021 | 0.67 | 0.45 |
| (6,8) | 0.25 | 0.36 | 0.017 | 0.44 | 0.25 |
| (7,3) | 0.40 | 0.40 | 0.021 | 0.62 | 0.40 |

(b)

**Figure 8.** (a) Visualization of *gdb19*, and (b) generated parameters for *gdb19-sdairp*.

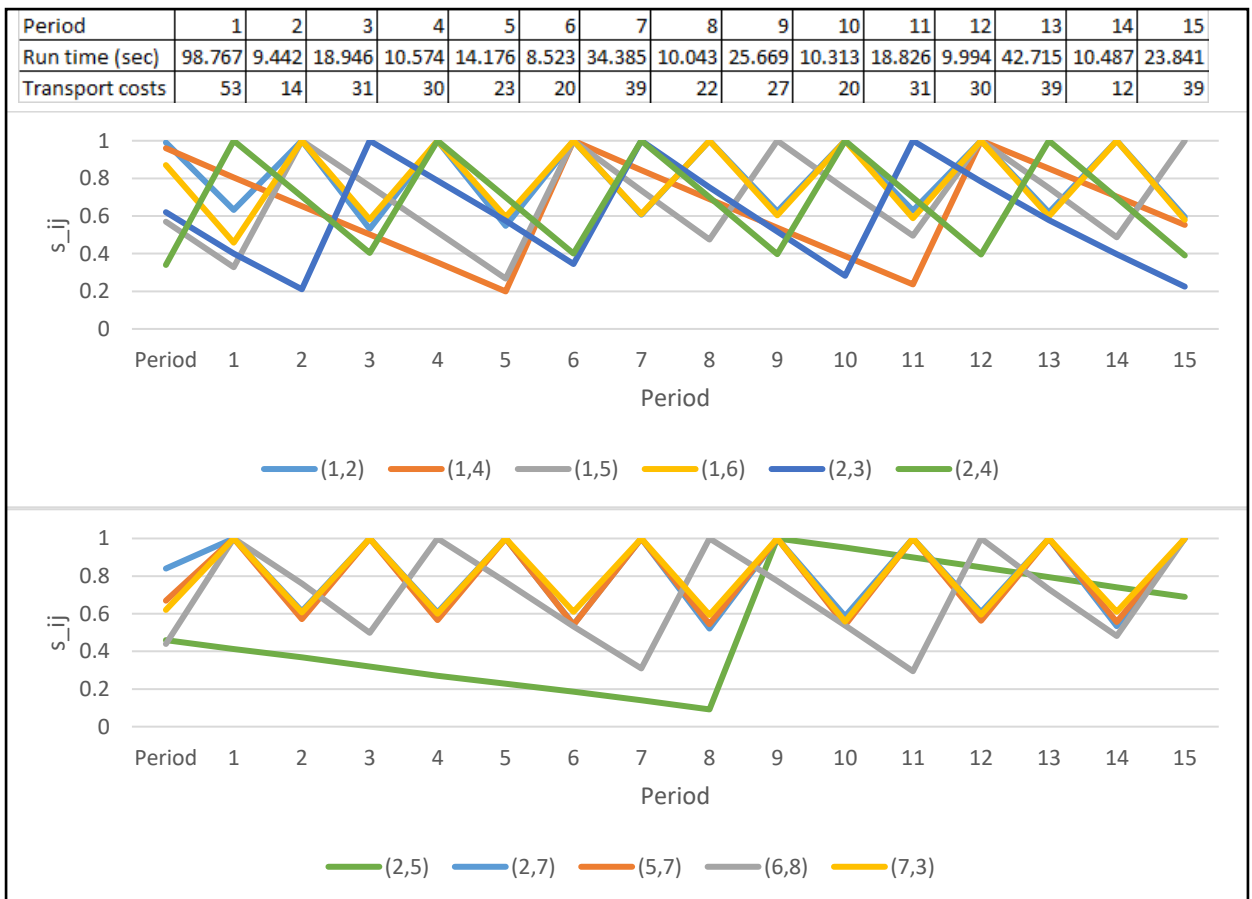

**Figure 9.** Summary of link inventories $s_{ij}$ over 15 periods as an outcome of proposed dynamic policy and algorithm.

Since the volatilities are small in this instance, the policy is essentially identical to the myopic policy. The run times for the algorithm in each period average to be 23.11 sec. Each time an inventory gets reset to 1, it indicates that a UAV was deployed to monitor the link. Some links, such as (2,5), take many periods to require attention, whereas others (e.g. (1,6)) need monitoring every two periods in this example. As a result of these differences, there is no convenient periodicity to set a static plan, and in this instance, coincidentally, no two periods had the same set of arcs monitored. This result proves why a naïve static multi-period deployment plan would not be effective.

## 6. Conclusion
### *6.1. Discussion for UAV city traffic monitoring*
Having established a new methodology for dynamically deploying UAVs, we discuss issues to apply this method for city monitoring. Mohammed et al. (2014) examined the challenges and opportunities of UAV deployment for smart cities, and concluded that these fall into multiple dimensions, including technology, governance, and management. For example, Foina et al. (2014) summarized regulations placed by the Federal Aviation Administration (FAA) on commercial UAV operations: for example, altitudes cannot exceed 500 ft. Several issues most relevant to the proposed methodology are discussed; these efforts are divided between scalability; stochastic variables; and drone technologies.

*6.1.1. Scale of problems*
While the method is demonstrated for a simple 5-node network over a range of different test settings for the purpose of replicability, city monitoring efforts will require deployment of dozens to hundreds of UAVs over a city. The scaling of the LSM algorithm itself generally is not an issue, as its computation time is based on $\{payoff\ computation\} \times \{number\ of\ time\ periods\} \times \{number\ of\ MC\ sampled\ paths\}$, and has been implemented in large scale in financial practice for option pricing. The scalability issue is largely in the $\{payoff\ computation\}$, which depends on running the internal selective VRP. As has been shown in Allahviranloo et al. (2014), the problem is NP-hard. As a result, testing this methodology on a much larger network would require implementing a heuristic algorithm on just the selective VRP solution, which would only add noise to the analysis of the effectiveness of the LSM approach.

In practice, the scale depends on the regulations in place limiting the feasible regions for monitoring, and criticality of corridors for monitoring, geographic scale of the region, etc. For example, a city may already have extensive sensors distributed over a city, and only require UAVs for certain critical corridors. This can be measured by calibrating the risk function parameter and stock out cost. Zhang et al. (2015) examined a different methodology of UAV for city monitoring with the Chicago network using only 4 UAVs. This suggests the complexity of the challenges in this issue.

Further complications can arise in a city space where the UAVs from multiple operators need to share (such as some for city monitoring while others are used for deliveries like Amazon's drones). This would have constraints on the feasible space for the problem.

A future study is needed to derive an optimal UAV fleet size problem for city monitoring depending on network structure, effectiveness of the monitoring information for incident management authorities, and volatility of the traffic conditions.

*6.1.2. Stochastic variables*

Traffic data was used as the real time information source for the SDAIRP. Other variables can also be incorporated. Short term weather conditions follow mean-reversion (e.g. Chow and Regan, 2011b), and may be used to restrict deployment (if weather is too unfavorable for safe deployment) or to encourage more deployment (icy conditions in post-blizzard settings). Other variables may include pedestrian counts, energy use in infrastructure, social media and location-based services (e.g. Twitter, Foursquare), etc. Regardless of the variable, it would need to be calibrated with the risk function.

*6.1.3. Drone technologies*

Ongoing efforts are advancing the technologies. The FAA has also deployed several test beds around the United States to evaluate the UAV technology, as shown in Figure 10. A number of commercial drones are now available on the market. According to Business Insider (2015), the market for commercial drones is expected to reach $3B by 2024. They cite several leaders in this market:
- senseFly (Switzerland)
- Aeryon (Canada)
- CybAero (Sweden)
- DJI (China)
- Gryphon (Korea)

Mohammed et al. (2014) highlight the Dragon Fly UAV used for surveillance. We are currently testing API functionality with a DJI Phantom 3 drone. These drones can be programmed to fly over GPS waypoints specified on mobile apps. Customizations can be made for designated launch sites.

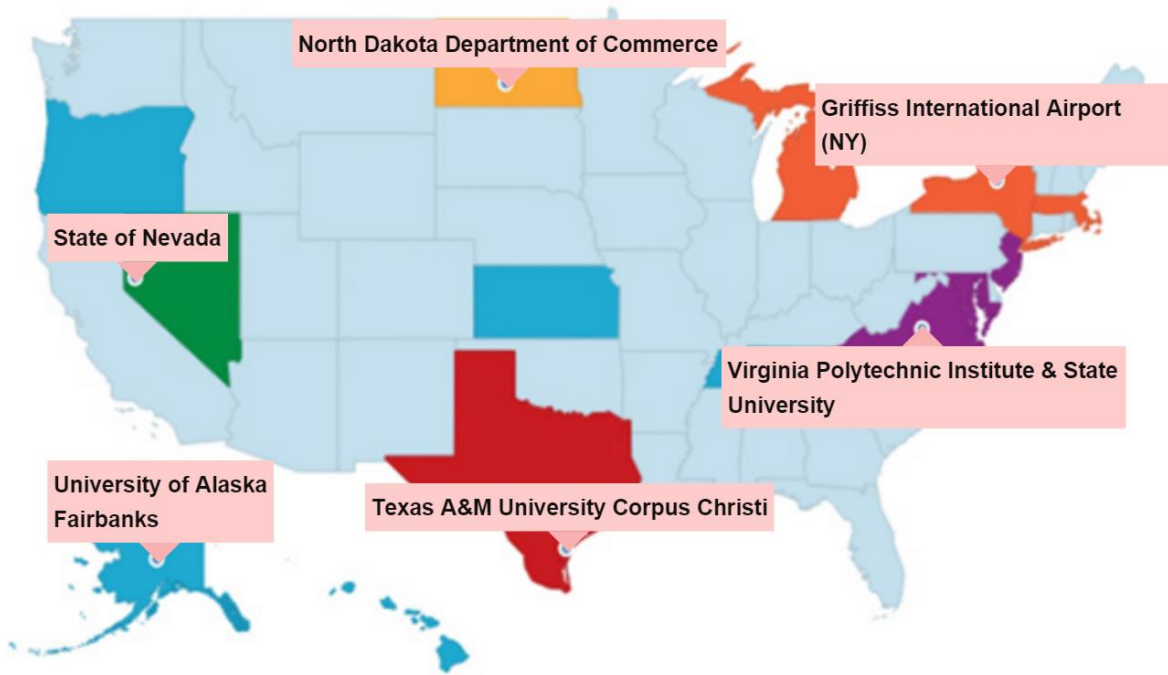

**Figure 10. FAA-approved test beds in the United States (source: FAA, 2015).**

*6.2. Summary and future directions*

We study the first deterministic arc-inventory routing problem and its stochastic dynamic policy, in an application to deploy mobile sensors like UAVs to monitor city traffic. The work is expected to contribute significantly to the growing literature on real time traffic monitoring, where uncertainty is highest and costliest.

More specifically, we address three research gaps in the literature. First, we synthesize and formulate the deterministic arc-inventory routing problem to handle multiple periods in such a way that can be modified for real time policies. The demand in the inventory routing problem is used to model the need for UAV monitoring based on data such as traffic conditions. The model is demonstrated with a replicable instance using a commercial solver.

We then consider the stochastic dynamic problem under real time data setting. Under mild assumptions about the stochastic demand being a stationary process, we formulate a model to approximate the optimal policy to the SDAIRP. The model features a newly defined risk function that is based on expected stock out costs in an order-up-to replenishment assumption. We further show that the recursive Bellman equation of the approximate policy can be formulated as a selective vehicle routing problem.

Lastly, we propose an approximate dynamic programming algorithm based on LSM, and modify the least squares estimation to be applied to the expected stock out cost upon the next replenishment as the "next time step". This new algorithm is tested on 30 simulated instances of real time trajectories over 5 time periods, resulting in solving 302,500 runs of the selective VRP to evaluate the proposed policy and algorithm. Computational results show that the algorithm on average outperforms the myopic policy by 23% to 28% depending on the parametric design. A second experiment is conducted using benchmark examples from the literature, where the instance *gdb19* is converted into a 15-period simulated scenario. A 5-period rolling horizon instance of the policy is run to illustrate why a naïve static multi-period deployment plan would not be effective, despite its current use in UAV deployment.

There are a number of steps to be taken in future research. Alternative solution methods based on policy iteration and value iteration should be compared with. The policy should be combined with a heuristic (e.g. parallel GA in Allahviranloo et al., 2014), a districting method, or continuous approximation for the underlying selective VRP for large scale deployment. While the multiple periods of recharging UAVs is considered in the deterministic formulations, it is not considered in the dynamic policy because we solve it using a backward dynamic programming procedure. Future research should seek to rectify this issue. Field tests for the different applications are also needed in traffic monitoring during major events to validate the performance of the deployment.

## Acknowledgements
This research was supported in part by the Canada Research Chairs program.

Wu, D., Arkhipov, D.I., Kim, M., Talcott, C.L., Regan, A.C., McCann, J.A., Venkatasubramanian, N., 2016. ADDSEN: adaptive data processing and dissemination for drone swarms in urban sensing. *IEEE Transactions on Computers*, in press, doi: 10.1109/TC.2016.2584061.

Yang, L., Chu, C.P., 2011. Stochastic model for traffic flow prediction and its validation. Proc. 90th Annual Meeting of the TRB, Washington DC.

Yazici, A., Kirlik, G., Parlaktuna, O., Sipahioglu, A., 2014. A dynamic path planning approach for multirobot sensor-based coverage considering energy constraints. *IEEE Trans. Cybernetics* 44(3), 305-314.

Zhang, J., Jia, L., Niu, S., Zhang, F., Tong, L., Zhou, X., 2015. A space-time network-based modeling framework for dynamic unmanned aerial vehicle routing in traffic incident monitoring applications. *Sensors* 15, 13874-13898.